\documentclass[conference]{IEEEtran}
\IEEEoverridecommandlockouts
\usepackage{cite}
\usepackage{amsmath,amssymb,amsfonts}
\usepackage{algorithmic}
\usepackage{graphicx}
\usepackage{textcomp}
\usepackage{xcolor}
\usepackage{comment}
\usepackage{caption}
\usepackage{epstopdf}
\usepackage{subcaption}

\usepackage[linesnumbered,ruled,vlined]{algorithm2e}
\SetKwInput{KwInput}{Input}                
\SetKwInput{KwOutput}{Output} 
\SetKwInput{KwInit}{Initialize}  
\SetKw{Break}{break}

\def\BibTeX{{\rm B\kern-.05em{\sc i\kern-.025em b}\kern-.08em
    T\kern-.1667em\lower.7ex\hbox{E}\kern-.125emX}}
\begin{document}

\title{Accelerated Proximal Iterative re-Weighted $\ell_1$  Alternating Minimization for Image Deblurring}

\author{\IEEEauthorblockN{Tarmizi Adam}
\IEEEauthorblockA{\textit{Faculty of Computing} \\
\textit{Universiti Teknologi Malaysia}\\
Jalan Iman, 81310, Skudai, Johor, Malaysia \\
tarmizi.adam@utm.my}
\and
\IEEEauthorblockN{Alexander Malyshev}
\IEEEauthorblockA{\textit{Department of Mathematics,} \\
	\textit{University of Bergen,}\\
	Postbox 7803, 5020 Bergen, Norway \\
	Alexander.Malyshev@uib.no}
\and
\IEEEauthorblockN{Mohd Fikree Hassan}
\IEEEauthorblockA{\textit{School of Information Technology,} \\
\textit{Monash University Malaysia,}\\
Jalan Lagoon Selatan, Bandar Sunway, \\
47500 Subang Jaya, Selangor, Malaysia \\
MohdFikree.Hassan@monash.edu}
\and
\IEEEauthorblockN{Nur Syarafina Mohamed}
\IEEEauthorblockA{\textit{Department of Mathematical Sciences, Faculty of Sciences,} \\
\textit{Universiti Teknologi Malaysia}\\
81310, Skudai, Johor, Malaysia \\
nursyarafina@utm.my }
\and
\IEEEauthorblockN{Md Sah Hj Salam}
\IEEEauthorblockA{\textit{Faculty of Computing} \\
	\textit{Universiti Teknologi Malaysia}\\
	Jalan Iman, 81310, Skudai, Johor, Malaysia \\
	sah@utm.my}
}

\maketitle

\begin{abstract}
The quadratic penalty alternating minimization (AM) method is widely used for solving the convex $\ell_1$ total variation (TV) image deblurring problem.  However, quadratic penalty AM for solving the nonconvex nonsmooth $\ell_p$, $0 < p < 1$ TV image deblurring problems is less studied. In this paper, we propose two algorithms, namely proximal iterative re-weighted $\ell_1$ AM (PIRL1-AM) and its accelerated version, accelerated proximal iterative re-weighted $\ell_1$ AM (APIRL1-AM) for solving the nonconvex nonsmooth $\ell_p$ TV image deblurring problem. The proposed algorithms are derived from the proximal iterative re-weighted $\ell_1$ (IRL1) algorithm and the proximal gradient algorithm. Numerical results show that PIRL1-AM is effective in retaining sharp edges in image deblurring while APIRL1-AM can further provide convergence speed up in terms of the number of algorithm iterations and computational time. 
\end{abstract}


\section{Introduction}
In this paper, we are interested in the image deblurring problem  obtained from  the following image degradation model
\begin{equation}
\label{eq:DegMod}
\mathbf{f} = \mathbf{Ku} + \mathbf{n},
\end{equation}
where $\mathbf{f}\in \mathbb{R}^n$ is the observed noisy and blurred image, $\mathbf{K} \in \mathbb{R}^{n \times n}$ is a blur kernel, $\mathbf{u} \in \mathbb{R}^n$ is the uncorrupted image to be estimated, and $\mathbf{n} \in \mathbb{R}^n$ is additive Gaussian noise. 

One way to solve (\ref{eq:DegMod}) for the image $\mathbf{u}$ is by minimizing the nonconvex nonsmooth composite optimization problem 
\begin{equation}
\label{eq:compMin}
 \underset{\mathbf{u} \in \mathbb{R}^n}{\text{min}}\,\frac{1}{2} \|  \mathbf{Ku} - \mathbf{f}  \|_2^2 + \mu  \| \nabla \mathbf{u}  \|_p^p,
\end{equation}
where $0 < p < 1$ and $\nabla \in \mathbb{R}^{n \times n}$ are the discrete difference operator \cite{lu2016implementation}. If $p = 1$, problem (\ref{eq:compMin}) results in the convex $\ell_1$-norm total variation (TV) image restoration \cite{rudin1992nonlinear}.



Operator splitting methods such as the alternating direction method of multipliers (ADMM) and the alternating minimization (AM) \cite{wang2008new,boyd2011distributed} for solving the convex $\ell_1$-norm TV problem produces the sublinear convergence rate of $\mathcal{O}\left( \frac{1}{k} \right)$ which is quite slow in practice\cite{he20121,beck2015convergence}. Therefore, this has motivated researchers to accelerate these methods to improved the convergence rate to $\mathcal{O}\left( \frac{1}{k^2} \right)$.  However, the majority of the work in this direction is focused on convex optimization.  

This paper, focuses on the model (\ref{eq:compMin}) when $0 < p < 1$ i.e., $\ell_p$ quasi-norm hence, the nonconvex nonsmooth $\ell_p$-norm TV. Our motivation mainly stems from the advantage of nonconvex nonsmooth penalties in restoring even sharper image quality compared to the $\ell_1$-norm TV \cite{zhang2020tv}. Furthermore, this paper is further motivated by applying the quadratic penalty AM which is an operator splitting method for solving nonconvex nonsmooth image deblurring problems (\ref{eq:compMin}).

By this, we propose a proximal iterative re-weighted $\ell_1$ alternating minimization (AM) algorithm along with its accelerated version.  The proposed algorithm uses the ideas of proximal operators and the iterative re-weighted $\ell_1$ (IRL1) method \cite{candes2008enhancing} in combination with the alternating minimization method \cite{wang2008new}. To accelerate the convergence of the proposed method, Nesterov acceleration is also used \cite{nesterov1983method, beck2009fast}. 


\section{Related Works}
The alternating minimization (AM) method for solving the $\ell_1$ TV image deblurring problem was initially proposed in \cite{wang2008new}. An accelerated version employing acceleration techniques from \cite{beck2009fast} was further investigated in \cite{xie2018new}. 


Specifically, it was shown that the $\ell_1$-norm sub minimization problem of the AM can be seen as a proximal gradient step. Hence, amenable to acceleration via the proximal gradient method \cite{beck2009fast}. However, the accelerated AM assumes the minimization problem to be convex.

The interesting link between AM and the accelerated proximal gradient method shown in \cite{xie2018new} suggests that AM may have links with nonconvex proximal gradient type methods. For example, the re-weighted $\ell_1$ method that was originally proposed in \cite{candes2008enhancing} for nonconvex $\ell_p$ sparse recovery problems along with its proximal version \cite{lu2014proximal,sun2017global}.

The proximal re-weighted $\ell_1$ algorithm with acceleration is relatively new and studied in \cite{wang2022,yu2019iteratively,phan2021accelerated} and has shown promising results in minimizing problems of the form (\ref{eq:compMin}). However, their applications have been mainly restricted to sparse and low-rank matrix recovery problems. Furthermore, their connections with operator splitting methods such as AM for image deblurring have not been to our knowledge explored.

\section{Problem Formulation and Algorithm}

In this paper, we focus on image deblurring by minimizing the nonconvex nonsmooth optimization problem (\ref{eq:compMin}).
 However, we restrict our results and discussion to a value of $p = 0.1$. By adding a quadratic penalty term to (\ref{eq:compMin}), we have
\begin{equation}
\label{eq:prob2}
\underset{\mathbf{u}}{\text{min}}\, \frac{1}{2} \| \mathbf{Ku} - \mathbf{f} \|_2^2 + \mu \| \mathbf{z} \|^p_p + \frac{\beta}{2}\| \mathbf{z} - \nabla \mathbf{u} \|_2^2.
\end{equation}
Indeed, problem (\ref{eq:prob2}) is equivalent to problem (\ref{eq:compMin}) with constraint $\mathbf{z} =  \nabla \mathbf{u}$ \cite{tan2014smoothing}. To minimize (\ref{eq:prob2}),  we can fix $\mathbf{u}$ with the current value and minimize with respect to $\mathbf{z}$ and vice-versa
\begin{equation}
\label{eq:AMitr}
\begin{cases}
\mathbf{z}_{k+1}= \underset{\mathbf{z}}{\operatorname{arg\,min}} \, \frac{\beta}{2}\| \mathbf{z} - \nabla \mathbf{u} \|_2^2 + \mu \| \mathbf{z} \|^p_p, \\
\mathbf{u}_{k+1}= \underset{\mathbf{u}}{\operatorname{arg\,min}} \, \frac{1}{2} \| \mathbf{Ku} - \mathbf{f} \|_2^2 + \frac{\beta}{2}\| \mathbf{z} - \nabla \mathbf{u} \|_2^2. \\
\end{cases}
\end{equation}
Quadratic penalty AM scheme (\ref{eq:AMitr}) is a classic method and commonly used in the image and signal processing literature \cite{tan2014smoothing,wang2008new}.

Note that the minimization problem concerning $\mathbf{z}$ in (\ref{eq:AMitr}) is the nonconvex nonsmooth $\ell_p$ minimization problem hence, amenable to the iterative re-weighted $\ell_1$ minimization algorithm \cite{candes2008enhancing}. If $p=1$ in (\ref{eq:AMitr}), this minimization problem is convex and solvable using soft thresholding \cite{wright2022optimization}. Consequently, it was shown that it is equivalent to the proximal gradient update \cite{xie2018new}
\begin{equation}
\begin{aligned}
\mathbf{z}_{k+1} = \underset{\mathbf{z} \in \mathbb{R}^n}{\text{argmin}}\, f\left( \mathbf{y}_k\right) + &\langle \nabla f\left( \mathbf{y}_k \right),\, \mathbf{z} - \mathbf{y}_k \rangle \\ 
&+ \frac{L}{2} \| \mathbf{z} - \mathbf{y}_k \|_2^2 + \mu \| \mathbf{z} \|_1,
\end{aligned}
\label{eq:pgm}
\end{equation}
where $\mathbf{y}_k = \nabla \mathbf{u}_k$, $L$ is the Lipschitz constant of $f\left( \mathbf{y}_k \right)=\frac{\beta}{2}\|\mathbf{z} - \mathbf{y}_k \|_2^2$, and $\nabla f\left( \mathbf{y}_k \right) = \beta \left(\mathbf{z} - \nabla \mathbf{u} \right)$. Due to this equivalence, the $\mathbf{z}$ sub problem can be accelerated by the fast iterative shrinkage and thresholding algorithm (FISTA) \cite{beck2009fast}.

\subsection{Proximal Iterative re-Weighted $\ell_1$ Alternating Minimization}
With previous foundations in place, coming back to the nonconvex nonsmooth $\ell_p$ minimization subproblem $\mathbf{z}$ in (\ref{eq:AMitr}), we have 
\begin{equation}
\mathbf{z}_{k+1}= \underset{\mathbf{z}}{\operatorname{arg\,min}} \, \frac{1}{2}\| \mathbf{z} - \nabla \mathbf{u}_k\|_2^2 + \frac{\mu}{\beta} \| \mathbf{z} \|^p_p,
\label{eq:subZa}
\end{equation}
after a simple re-arrangement of $\beta$. The IRL1 approximately solves (\ref{eq:subZa}) by \cite{candes2008enhancing,wang2021relating} 
\begin{equation}
\mathbf{z}_{k+1} =\underset{\mathbf{z}\in \mathbb{R}^n}{\operatorname{arg\,min}} \, \frac{1}{2}\| \mathbf{z} - \nabla \mathbf{u}_k \|_2^2 +\frac{\mu}{\beta} \sum_i w^i | z^i |,
\label{eq:sumIrl1}
\end{equation}
where $w^i=p\left( |z^i_k| + \epsilon \right)^{p-1}$, and $z^i$ are the weights and entries of vector $\mathbf{z}$ respectively. Note that in (\ref{eq:sumIrl1}) the nonconvex nonsmooth $\ell_p$ minimization problem is approximated into a convex weighted $\ell_1$-norm minimization hence, IRL1 is a convex relaxation method for nonconvex optimization problems.

By introducing a diagonal weight matrix $\mathbf{W}_k = \text{diag}\left( w^i_k, \cdots w_k^n \right)$, IRL1 problem (\ref{eq:sumIrl1}) can be written as
\begin{equation}
\mathbf{z}_{k+1} =\underset{\mathbf{z}\in \mathbb{R}^n}{\operatorname{arg\,min}} \, \frac{1}{2}\| \mathbf{z} - \nabla \mathbf{u}_k \|_2^2 + \frac{\mu}{\beta}\| \mathbf{W}_k \mathbf{z} \|_1.
\label{eq:matIrl1}
\end{equation}
Problem (\ref{eq:matIrl1}) can be interpreted as solving a proximal linearization of the term $f\left( \mathbf{z} \right)=\frac{1}{2}\| \mathbf{z} - \bar{\mathbf{z}_k} \|_2^2 $ at $\bar{\mathbf{z}_k}= \nabla \mathbf{u}_k$ i.e.,  
\begin{equation}
\begin{aligned}
\mathbf{z}_{k+1} = \underset{\mathbf{z} \in \mathbb{R}^n}{\operatorname{arg\,min}}\, f\left( \bar{\mathbf{z}_k}\right) + &\langle \nabla f\left( \bar{\mathbf{z}_k} \right),\, \mathbf{z} - \bar{\mathbf{z}_k} \rangle \\
&+ \frac{L}{2} \| \mathbf{z} - \bar{\mathbf{z}_k} \|_2^2 + \frac{\mu}{\beta} \| \mathbf{W}_k \mathbf{z} \|_1,
\end{aligned}
\label{eq:NvxPgm}
\end{equation}
hence, by ignoring constant terms can be written as \cite{beck2017first,beck2009fast}
\begin{equation}
\mathbf{z}_{k+1} =\underset{\mathbf{z}\in \mathbb{R}^n}{\operatorname{arg\,min}} \, \frac{1}{2}\| \mathbf{z} -  \mathbf{v}_k \|_2^2 + \lambda \| \mathbf{W}_k \mathbf{z} \|_1,
\end{equation}
with $\mathbf{v}_k = \bar{\mathbf{z}_k} - \frac{1}{L}\nabla f\left(\bar{\mathbf{z}_k} \right) $ and $\lambda = \frac{\mu}{L\beta}$. Furthermore, by combining the weight entries with $\lambda$ in the weight matrix $\mathbf{W}_k$ we finally arrive at
\begin{equation}
\mathbf{z}_{k+1} =\underset{\mathbf{z}\in \mathbb{R}^n}{\operatorname{arg\,min}} \, \frac{1}{2}\| \mathbf{z} -  \mathbf{v}_k \|_2^2 + \frac{p\lambda}{\left( |\bar{z}^i_k| + \epsilon \right)^{1-p}} \| \mathbf{z} \|_1,
\label{eq:FinUpd}
\end{equation}
where $\bar{z}^i_k$ is the $i^\text{th}$ entry of $\bar{\mathbf{z}}_k$ at the  $k^\text{th}$ iteration and $\epsilon$ a very small number to avoid division by zero. Equation (\ref{eq:FinUpd}) can be solved in a closed form via the soft-thresholding operation as
\begin{equation}
\mathbf{z}_{k+1} = \text{sgn}\left(\bar{z}_k^i  \right) \text{max} \left( 0,\, |\bar{z}_k^i| - \frac{p\lambda}{\left( |z^i_k| + \epsilon \right)^{1-p}} \right).
\label{eq:lpShrink}
\end{equation}
From (\ref{eq:lpShrink}), it can be seen that solving the original nonconvex subproblem (\ref{eq:subZa}) boils down to a series of proximal weighted $\ell_1$ minimization problem. 

The idea of proximal re-weighted $\ell_1$ minimization has been proposed in \cite{lu2014proximal} and its convergence behavior analyzed in \cite{sun2017global}. However, the applications are mainly confined to sparse signal recovery and its use as a sub-minimization problem in the AM scheme (\ref{eq:AMitr}) for image deblurring to our knowledge has not been studied.

Next, the sub-minimization problem for $\mathbf{u}$
\begin{equation}
\mathbf{u}_{k+1}= \underset{\mathbf{u}}{\operatorname{arg\,min}} \, \frac{1}{2} \| \mathbf{Ku} - \mathbf{f} \|_2^2 + \frac{\beta}{2}\| \mathbf{z}_{k+1} - \nabla \mathbf{u} \|_2^2,
\end{equation}
is a convex quadratic problem that can be solved by solving the following linear system
\begin{equation}
\mathbf{u}_{k+1} = \left( \mathbf{K}^\top \mathbf{K} + \beta \nabla^\top \nabla \right)^{-1} \left( \mathbf{K}^\top \mathbf{f} + \beta \nabla^\top \mathbf{z}_{k+1} \right).
\label{eq:linSyst}
\end{equation}
Taking into account equations (\ref{eq:lpShrink}) and (\ref{eq:linSyst}), the complete listing for the proximal iterative re-weighted $\ell_1$ alternating minimization (PIRL1-AM) is listed as Algorithm 
\ref{al:am1}.
\begin{algorithm}
\caption{Proximal iterative re-weighted $\ell_1$ alternating minimization (PIRL1-AM)}
\label{al:am1}

\KwInit{$\mathbf{u}_0$, $\mathbf{z}_0$, $\mu>0$, $\beta > 0$, $p = 0.1 $, $L = 1$, and $k = 0$}

\While{not converged}
{

$\mathbf{z}_{k+1} = \text{sgn}\left(\bar{z}_k^i  \right) \text{max} \left( 0,\, |\bar{z}_k^i| - \frac{p\lambda}{\left( |z^i_k| + \epsilon \right)^{1-p}} \right)$, \\

  $\mathbf{u}_{k+1} = \left( \mathbf{K}^\top \mathbf{K} + \beta \nabla^\top \nabla \right)^{-1} \left( \mathbf{K}^\top \mathbf{f} + \beta \nabla^\top \mathbf{z}_{k+1} \right)$,\\
$k= k+1$
}
\end{algorithm}

\subsection{Accelerated Proximal Iterative re-Weighted $\ell_1$ Alternating Minimization}
To accelerate the proximal iterative re-weighted $\ell_1$ AM for the nonconvex nonsmooth $\ell_p$ TV image deblurring problem, accelerated techniques from \cite{nesterov1983method, beck2009fast} can be used.
Consider the sub-minimization problem (\ref{eq:FinUpd}),
due to the convexity of (\ref{eq:FinUpd}) and its equivalence to minimizing a proximal linearization of $\frac{1}{2}\| \mathbf{z} -  \mathbf{v}_k \|_2^2 $ as discussed earlier, we have
\begin{equation}
\begin{aligned} 
\label{eq:ProxACCIR}
\mathbf{z}_{k+1} = \underset{\mathbf{z} \in \mathbb{R}^n}{\text{argmin}}\, f\left( \mathbf{v}_k\right) + &\langle \nabla f\left( \mathbf{v}_k \right),\, \mathbf{z} - \mathbf{v}_k \rangle \\ 
&+\frac{L}{2} \| \mathbf{z} - \mathbf{v}_k \|_2^2 + \tau \|  \mathbf{z} \|_1,
\end{aligned}
\end{equation}
where $\tau =  \frac{p\lambda}{\left( |\bar{z}^i_k| + \epsilon \right)^{1-p}} $. Then, the acceleration strategies along the lines of \cite{beck2009fast,nesterov1983method} can be employed, which involves the following iterative scheme
\begin{equation}
\begin{aligned} 
\begin{cases}
\mathbf{z}_{k+1} = \underset{\mathbf{z} \in \mathbb{R}^n}{\text{argmin}}\, f\left( \mathbf{y}_k\right) + \langle \nabla f\left( \mathbf{y}_k \right),\, \mathbf{z} - \mathbf{y}_k \rangle \\
\qquad \qquad \qquad \qquad \qquad +\frac{L}{2} \| \mathbf{z} - \mathbf{y}_k \|_2^2 + \tau \|  \mathbf{z} \|_1,\\
 t_k = \frac{k-1}{k + 2}, \\
\mathbf{y}_{k+1} = \mathbf{z}_{k+1} + t_k \left( \mathbf{z}_{k+1} - \mathbf{z}_k \right).\\
\end{cases}
\end{aligned}
\label{eq:itrScmACCIRL1}
\end{equation}
The $\mathbf{z}$ minimization in (\ref{eq:itrScmACCIRL1}) as disccused previously has a closed form solution of
\begin{equation}
\mathbf{z}_{k+1} = \text{sgn}\left(\bar{z}_k^i  \right) \text{max} \left( 0,\, |\bar{z}_k^i| - \tau \right).
\end{equation}
In the iterative scheme (\ref{eq:itrScmACCIRL1}), the scalar $t_k$ is known as the Nesterov momentum coefficient and changes in each iteration $k$. The step $\mathbf{y}_k$ is called the extrapolation step. Also, recall that $\mathbf{z}_k = \nabla\mathbf{u}_k$. Acceleration techniques of this form are shown to match the theoretical lower bound of $\mathcal{O}\left( \frac{1}{k^2} \right)$ for smooth first-order convex optimization. 

Having shown the equivalence between (\ref{eq:FinUpd}) and the proximal linearization step (\ref{eq:ProxACCIR}) along with its acceleration, the iterative scheme of the accelerated proximal iterative re-weighted $\ell_1$ AM (APIRL1-AM) taking account the sub minimization problem for $\mathbf{u}$ is as follows
\begin{equation}
\begin{cases}
\mathbf{z}_{k+1} = \text{sgn}\left(\bar{z}_k^i  \right) \text{max} \left( 0,\, |\bar{z}_k^i| - \tau \right),\\
t_k = \frac{k-1}{k + 2}, \\
\mathbf{y}_{k+1} = \mathbf{z}_{k+1} + t_k \left( \mathbf{z}_{k+1} - \mathbf{z}_k \right),\\
\mathbf{u}_{k+1} = \left( \mathbf{K}^\top \mathbf{K} + \beta \nabla^\top \nabla \right)^{-1} \left( \mathbf{K}^\top \mathbf{f} + \beta \nabla^\top \mathbf{z}_{k+1} \right).\\
\end{cases}
\label{eq:PropAcc}
\end{equation}
The complete algorithm is listed as Algorithm \ref{al:amProp}. 

\begin{algorithm}
\caption{Accelerated proximal iterative re-weighted $\ell_1$ alternating minimization (APIRL1-AM)}
\label{al:amProp}

\KwInit{$\mathbf{u}_0$, $\mathbf{z}_0$, $\mu>0$, $\beta > 0$, $p = 0.1 $, $L = 1$, and $k = 0$}

\While{not converged}
{

$\mathbf{z}_{k+1} = \text{sgn}\left(\bar{z}_k^i  \right) \text{max} \left( 0,\, |\bar{z}_k^i| - \tau \right)$\\
$t_k = \frac{k-1}{k + 2}$ \\
$\mathbf{y}_{k+1} = \mathbf{z}_{k+1} + t_k \left( \mathbf{z}_{k+1} - \mathbf{z}_k \right)$\\
$\mathbf{u}_{k+1} = \left( \mathbf{K}^\top \mathbf{K} + \beta \nabla^\top \nabla \right)^{-1} \left( \mathbf{K}^\top \mathbf{f} + \beta \nabla^\top \mathbf{z}_{k+1} \right)$.\\
$k= k+1$
}
\end{algorithm}

In Algorithms \ref{al:am1} and \ref{al:amProp}, the Lipschitz constant $L$ is fixed to 1. In real applications, the value of $L$ is usually unknown. The values of $\mu$ and $\beta$ are used for computing $\lambda = \frac{\mu}{L\beta}$.
Additionally, the linear system for solving $\mathbf{u}$ can be solved very fast using the two-dimensional fast Fourier transforms (FFT) with complexity $\mathcal{O}\left(n^n \log n \right)$ \cite{wang2008new} while the $\mathbf{z}$ problem is only of linear complexity $\mathcal{O}\left( n \right)$.

\section{Results}
In this section, we apply the two proposed algorithms on the nonconvex nonsmooth image deblurring problem (\ref{eq:compMin}) and discuss the results of the proposed algorithms without and with acceleration.
\subsection{Experimental Setup}
For the experiments, the blurred and noisy images are obtained by the model (\ref{eq:DegMod}). The blur kernel used is a Gaussian blur of size $17\times 17$ pixel with $\sigma = 7$. Two blur levels were used namely, blur signal-to-noise ratio (BSNR) 30 and 20 \cite{banham1997digital}. The values for $\beta$ are $\beta = 0.01$  and $\beta = 0.009$ for BSNR 20 and 30 respectively. The deblurring on each image was done 10 times and the average was taken.

All images are of size $512\times 512$ pixels. Experiments were conducted using MATLAB\footnote{Codes at https://github.com/tarmiziAdam2005/PIRL1-AM} on an Intel Core i3-10105 CPU operating at 3.70 GHz with a memory of 4 Gb.

\begin{table*}[]
\caption{Image deblurring Results. Value of $\mu= 30,\, 60$ for Peppers and Cameraman respectively for BSNR 30. For BSNR 20, Value of $\mu= 100$ for both Peppers and Cameraman}
\centering
\begin{tabular}{cccccccccccc}
\hline
 &  &  &  &APIRL1-AM   & &  &&  & PIRL1-AM    \\ 
\hline
 &BSNR lvl  &  &PSNR  &SSIM  &Itr  &T &  &PSNR  &SSIM  &Itr  &T   \\
\hline
 &  &\textit{Peppers}  &27.39  &0.759   &269 &11.20 &  &27.39 &0.759 &364 &14.71  \\
 &30  &\textit{Cameraman}  &26.03 &0.763  &228 &9.50 &  &26.04 &0.762 &324 &13.12  \\
 &  &  &  &  &  &  &  & & &   \\
 &  &\textit{Peppers}  &25.52  &0.625   &141 &6.03 &  &25.52 &0.625  &186 &7.63  \\
 &20  &\textit{Cameraman}  &24.49  &0.570   &225 &9.64 &  &24.49 &0.570 &244 &10.02  \\
\hline
\end{tabular}
\label{tab: gausRes}
\end{table*}

\begin{figure}
 \centering
     \begin{subfigure}[b]{0.5\textwidth}
        \centering
        \includegraphics[height=5.1cm, width=8.3cm]{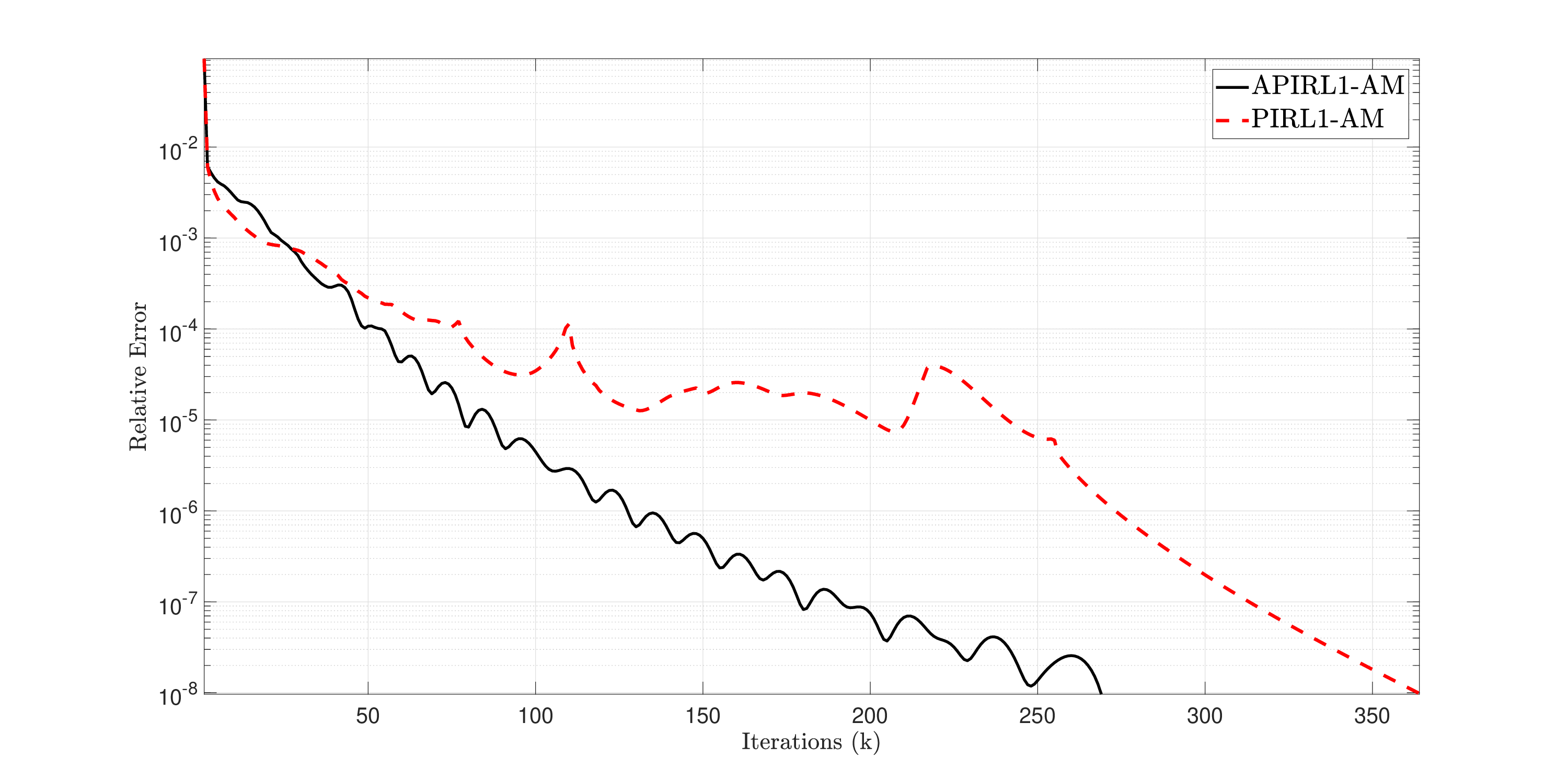}
        \caption{}
    \end{subfigure} 
    \begin{subfigure}[b]{0.5\textwidth}
        \centering
        \includegraphics[height=5.1cm, width=8.3cm]{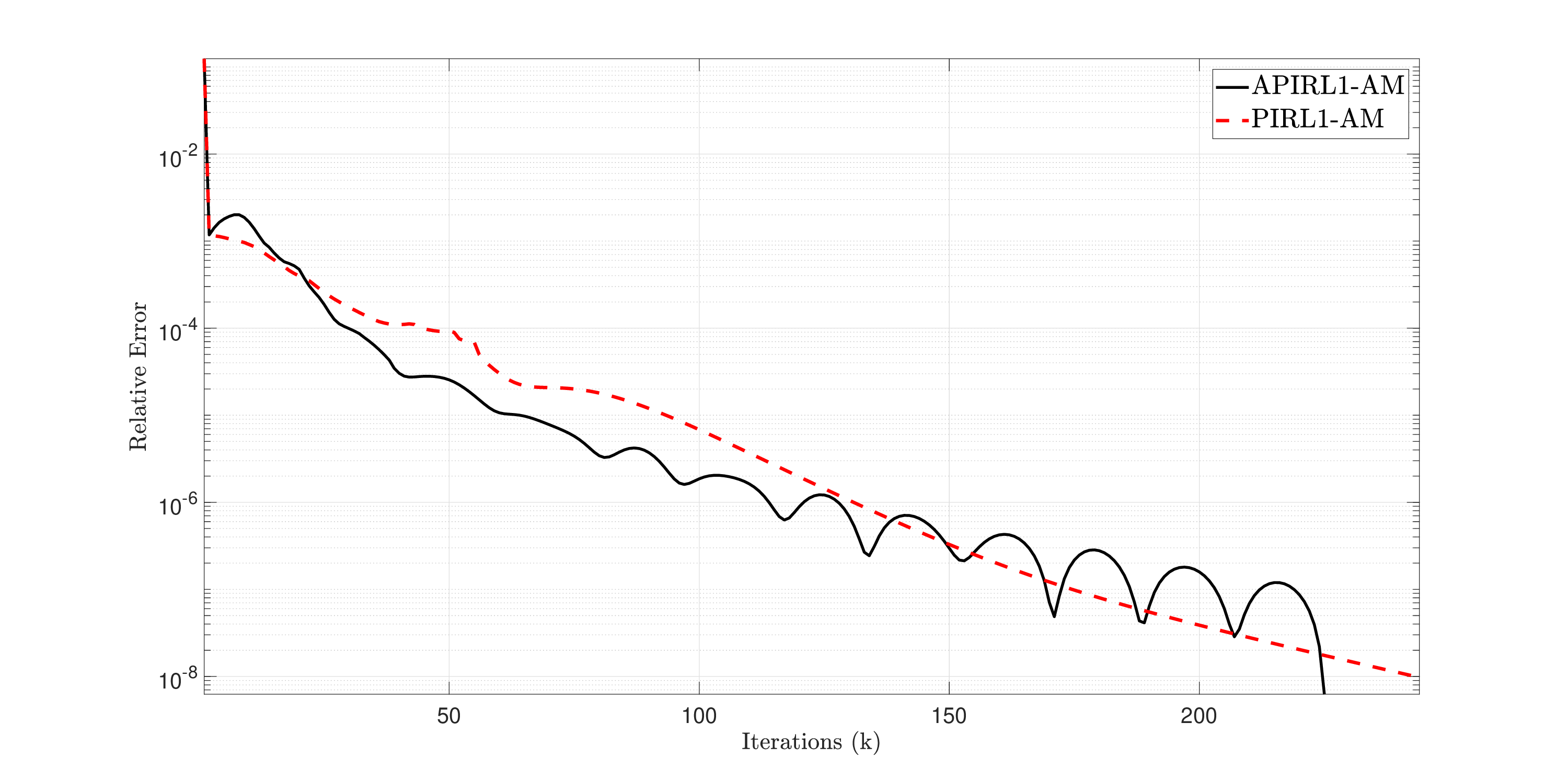}
         \caption{}
    \end{subfigure} 
\caption{Convergence of the relative errors $\frac{\| \mathbf{u}_{k} - \mathbf{u}_{k+1} \|_2}{\| \mathbf{u}_k \|_2} \leq 1\times 10^{-8}$, produced by the iterates of APIRL1-AM  and PIRL1-AM  for imagem (a) \textit{Peppers} (BSNR $=30$ ) and (b) \textit{Cameraman} (BSNR $=20$).}
\label{fig:convergencePep30}
\end{figure}

\begin{figure}
 \centering
    \begin{subfigure}[b]{0.14\textwidth}
        \centering
        \includegraphics[height=2.3cm, width=2.5cm]{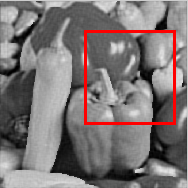}
        \caption{}
    \end{subfigure}
    \begin{subfigure}[b]{0.14\textwidth}
        \centering
        \includegraphics[height=2.3cm, width=2.5cm]{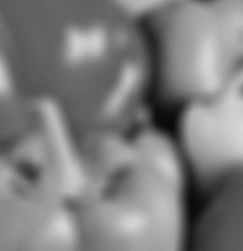}
        \caption{}
    \end{subfigure} 
\\[1.4ex]
\begin{subfigure}[b]{0.14\textwidth}
        \centering
        \includegraphics[height=2.3cm, width=2.5cm]{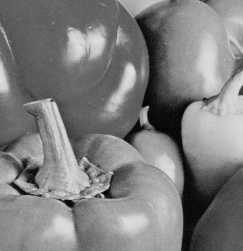}
        \caption{}
    \end{subfigure}
\begin{subfigure}[b]{0.14\textwidth}
        \centering
        \includegraphics[height=2.3cm, width=2.5cm]{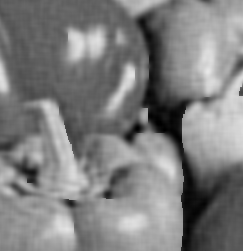}
        \caption{}
    \end{subfigure}
 \caption{Deblurred image \textit{Peppers} BSNR$ = 30$. (a) Original. (b) Blurred (PSNR = 22.97, SSIM = 0.659). (c) Zoomed. (d) APIRL1-AM deblurred.} 
  \label{fig:test2} 
\end{figure}

\begin{figure}
 \centering
    \begin{subfigure}[b]{0.14\textwidth}
        \centering
        \includegraphics[height=2.3cm, width=2.5cm]{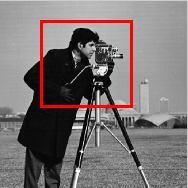}
        \caption{}
    \end{subfigure}
    \begin{subfigure}[b]{0.14\textwidth}
        \centering
        \includegraphics[height=2.3cm, width=2.5cm]{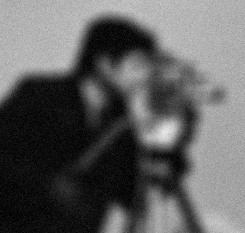}
        \caption{}
    \end{subfigure} 
\\[1.4ex]
\begin{subfigure}[b]{0.14\textwidth}
        \centering
        \includegraphics[height=2.3cm, width=2.5cm]{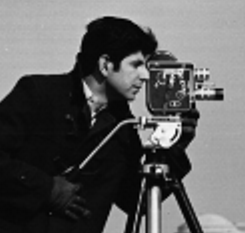}
        \caption{}
    \end{subfigure}
\begin{subfigure}[b]{0.14\textwidth}
        \centering
        \includegraphics[height=2.3cm, width=2.5cm]{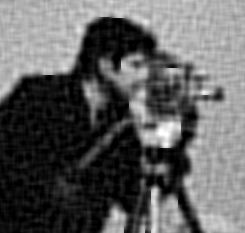}
        \caption{}
    \end{subfigure}
 \caption{Deblurred image \textit{Cameraman} BSNR$ = 20$. (a) Original. (b) Blurred (PSNR = 21.65, SSIM = 0.499). (c) Zoomed. (d) APIRL1-AM deblurred.} 
  \label{fig:test1} 
\end{figure}

\subsection{Results and discussion}
Table \ref{tab: gausRes} shows the results of the two proposed algorithms. In terms of image quality of PSNR and SSIM \cite{wang2004image}, the algorithms gave almost similar results. Some\footnote{Only deblurred images of APIRL1-AM are shown due to both algorithms giving very similar PSNR and SSIM values.} of the deblurred images are shown in Figures \ref{fig:test2} and \ref{fig:test1}.
For BSNR 30 (less noise corruption), Figure \ref{fig:test2} shows the ability of the methods in preserving sharp images.   

However, there is a difference in terms of the number of iterations to converge to the required relative error value. For both BSNR levels tested,  the acceleration technique used in APIRL1-AM manages to decrease the number of iterations to converge. Additionally, the time to converge also improves for APIRL1-AM compared to PIRL1-AM. 

Figure \ref{fig:convergencePep30} compares the convergence between APIRL1-AM and PIRL1-AM. It can be seen that APIRL1-AM converges faster than PIRL1-AM and exhibits acceleration ripples akin to first-order accelerated techniques \cite{o2015adaptive}. Taking into account similar results in terms of PSNR and SSIM of the two algorithms, acceleration gives an additional advantage for arriving at similar deblurring quality at a lower number of iterations and CPU time.

\section{Conclusion}
In this paper, two algorithms PIRL1-AM and APIRL1-AM were proposed for nonconvex nonsmooth $\ell_p$ TV image deblurring. The algorithms were derived by showing the links between the proximal gradient method and the proximal iterative re-weighted $\ell_1$. Both algorithms were able to retain sharp images for image deblurring. For algorithm convergence, APIRL1-AM exhibits the optimal $\mathcal{O}\left( \frac{1}{k^2}\right)$ rate of convergence and improves the CPU time and the number of iterations to converge. Future works include establishing the convergence rate of both algorithms and applying them to different problems.

\bibliographystyle{IEEEtran}
\bibliography{ICIP2023}

\end{document}